\documentclass[12pt]{amsart}
\usepackage{amsfonts}
\usepackage{latexsym}
\usepackage{amscd}
\usepackage{amsmath}
\usepackage{amssymb}
\usepackage{graphicx}
\usepackage{xypic}
\usepackage[mathscr]{eucal}
\usepackage[dvips]{color}

\newtheorem{de}{Definition}
\newtheorem{pro}{Proposition}
\newtheorem{cor}{Corollary}
\newtheorem{teo}{Theorem}
\newtheorem{rem}{Remark}
\newtheorem{lem}{Lemma}
\newtheorem{exa}{Example}
\newtheorem{alg}{Algorithm}

\newcommand{\bd}{\begin{de} \rm }
\newcommand{\ed}{ \end{de}}

\newcommand{\bp}{\begin{pro}}
\newcommand{\ep}{ \end{pro}}

\newcommand{\bc}{\begin{cor}}
\newcommand{\ec}{ \end{cor}}

\newcommand{\bt}{\begin{teo}}
\newcommand{\et}{ \end{teo}}

\newcommand{\br}{\begin{rem} \rm}
\newcommand{\er}{ \end{rem}}

\newcommand{\bl}{\begin{lem}}
\newcommand{\el}{ \end{lem}}

\newcommand{\be}{\begin{exa} \rm }
\newcommand{\ee}{ \end{exa}}

\newcommand{\ba}{\begin{alg}}
\newcommand{\ea}{ \end{alg}}

\newcommand{\bi}{\begin{itemize}}
\newcommand{\ei}{ \end{itemize}}

\newcommand{\co}{{\mathcal O}}

\newcommand{\gp}{\mathbb{P}}

\newcommand{\gr}{\mathbb{R}}

\newcommand{\findemo}{$\ \ \square$}

\newcommand{\pic}{{\rm Pic}}

\title{Fibers of pencils of curves on smooth surfaces}
\author{Francisco Monserrat} \curraddr{Instituto Universitario de
Matem\'atica Pura y Aplicada, Universidad Polit\'ecnica de Valencia,
Valencia  (Spain)} \email{framonde@mat.upv.es}
\date{}
\thanks{Supported by Spain Ministry of Education
 MTM2007-64704, JCyL VA025A07 and Bancaixa P1-1A2005-08}

\begin{document}

\maketitle

\begin{abstract}
Let $X$ be a smooth projective surface such that linear and
numerical equivalence of divisors on $X$ coincide and let
$\sigma\subseteq |D|$ be a linear pencil on $X$ with integral
general fibers. A fiber of $\sigma$ will be called special if
either it is not integral or it has non-generic multiplicity at
some of the base points (including the infinitely near ones) of
the pencil. In this note we provide an algorithm to compute the
integral components of the special fibers of $\sigma$.
\end{abstract}

\section{Introduction}

In this note we consider a (linear) pencil of curves without fixed
components on a smooth projective surface $X$ over an algebraically
closed field $k$ of arbitrary characteristic. Our goal is to provide
an algorithm, that uses infinitely near points, for computing the
integral (reduced and irreducible) components of the {\it special
fibers} of the pencil, assuming that it has integral general fibers
and that linear and numerical equivalence of divisors on $X$
coincide (condition which is satisfied, for instance, by projective
smooth rational surfaces). The concept of {\it special} fiber we
shall use will be defined later (in particular, the reducible fibers are
special).

There are several results in the literature related with special fibers of a pencil.
In \cite{vistoli, bodin, pereira, yuz}, bounds either on the
number of reducible fibers or on the number of reducible fibers of a particular
type are provided. The articles \cite{lib} and \cite{falk} are
the starting point of a new motivation for the interest
of completely reducible fibers of a pencil of projective plane curves (that is, unions of lines not necesarily
reduced) due to their relationship with the theory of line arrangements.
In \cite{dimca}, certain special fibers of pencils appear in the study of characteristic
varieties of local systems on a plane curve arrangement complement.

We consider an effective divisor $D$
on $X$ and a linear system $\sigma \subseteq |D|$ without fixed
components and with projective dimension 1 (a {\it pencil} on $X$,
in the sequel). It
corresponds to the projectivization of the sub-vector space
$V_{\sigma}$ of $H^0(X,\co_X(D))$ given by $\{s \in
H^0(X,\co_X(D)) \mid (s)_0\in \sigma\}\cup \{0\},$ $(s)_0$
denoting the divisor of zeros of the section $s$ (see
\cite[II.7]{har}). A basis of $V_{\sigma}$ provides a rational map $f_{\sigma}: X
\cdots \rightarrow \gp_k^1$, which is independent from the basis
up to composition with an automorphism of $\gp_k^1$. The closures
of the fibers of $f_{\sigma}$ are exactly the curves of the pencil
$\sigma$.  Moreover, there exists a minimal composition
of point blowing-ups $\pi_{\sigma}: Z_{\sigma} \longrightarrow X$
eliminating the indeterminacies of
the rational map $f_{\sigma}$, that is,
the map $h_{\sigma}:= f_{\sigma} \circ
\pi_{\sigma}:Z_{\sigma}\rightarrow \gp_k^1$ is a morphism.
We denote by $BP(\sigma)$ to the set of centers of the
blowing-ups used to get it; $Z_{\sigma}$ is the so-called {\it sky} of $BP(\sigma)$. The morphism $\pi_{\sigma}$ and the set $BP(\sigma)$ are essentially unique, since we identify those corresponding to $X$-isomorphic skies.

For each point $p \in BP(\sigma)$, set $E^*_p$ the total transform
on $Z_\sigma$ of the exceptional divisor created by the blowing-up
at $p$. Now, assume that $C$ is a curve on $X$. Then,
$$\pi_{\sigma}^*C=\tilde{C}+\sum_{p\in BP(\sigma)} m_p(C) E_p^*,$$
where $\tilde{C}$ denotes the  strict transform of $C$ on
$Z_{\sigma}$ and $m_p(C)$ is the multiplicity at $p$ of the strict
transform of $C$ on the surface to which $p$ belongs. The
successive strict transforms of the fibers of the pencil $\sigma$,
except a finite number, have the same multiplicity at every point
$p\in BP(\sigma)$, that we shall denote by $m_p(\sigma)$.  The
strict transforms on $Z_{\sigma}$ of these (infinitely many)
fibers belong to the complete linear system $|G_{\sigma}|$, where
$G_{\sigma}:= \pi_{\sigma}^*D-\sum_{p\in BP(\sigma)} m_p(\sigma)
E_p^*$. If the pencil $\sigma$ is {\it irreducible} (that is, with reduced
and irreducible general fibers), we shall say that a fiber $C$ of
$\sigma$ is {\it special} if either it is not integral or there
exists $p\in BP(\sigma)$ such that $m_p(C)\not= m_p(\sigma)$.

\section{Computation of the components of the special fibers}

Consider, as above, a smooth projective surface $X$ over $k$ such
that linear and numerical equivalence of divisors coincide. Fix a
closed immersion $i: X \rightarrow \gp^s_k$. For each curve $C$ on
$X$ (resp., its linear equivalence class ${\mathcal L}$), $\deg C$
(resp., $\deg {\mathcal L}$) will denote the degree of $C$ (resp.,
${\mathcal L}$) with respect to the immersion $i$ (that is, the
intersection product $i^*\co_{\gp^s_k}(1)\cdot C$). For each divisor
$D$ on a surface $Z$, $[D]$ will denote its class in $\pic(Z)$ (and
also in $\pic(Z)\otimes \gr$ via the natural inclusion map). Next,
we shall prove a statement that will support our method to compute
the integral components of the special fibers of an irreducible
pencil.

\bp \label{gordo} Let $\sigma\subseteq |D|$ be a pencil and $C$ an
integral curve, both on $X$. Then, $C$ is a component of a fiber
of $\sigma$ if, and only if, $G_{\sigma}\cdot \tilde{C}=0$.
Moreover, in this case, $\tilde{C}^2\leq 0$. \ep

\noindent {\it Proof.} If $C$ is a component of a fiber of
$\sigma$, then $\tilde{C}$ does not meet a general fiber of
$h_{\sigma}$ and so $G_{\sigma}\cdot \tilde{C}=0$. Conversely, if
$C$ is not included in a fiber of $\sigma$, then
$h_{\sigma}(\tilde{C})=\gp^1$ and $\tilde{C}$ meets properly each
fiber; therefore $G_{\sigma}\cdot \tilde{C}>0$. For the last
assertion observe that, as a consequence of the Hodge Index
Theorem \cite[IV.1.9]{har}, the set of elements $x$ of the real
vector space $\pic(Z_{\sigma})\otimes \gr$ such that $x\cdot x\geq
0$ and $[H]\cdot x\geq 0$, where $H$ is an ample divisor, (considering the bilinear pairing
induced by the intersection product of divisors) is the half-cone
over an Euclidean ball of dimension $\rho(Z_{\sigma})-1$,
$\rho(Z_{\sigma})$ being the Picard number of $Z_{\sigma}$; the
orthogonal hyperplane to $[G_{\sigma}]$ is tangent to that
half-cone, since $G_{\sigma}^2=0$. Taking these facts into
account, the second assertion follows from the first one and the
strict
convexity of the Euclidean ball.\findemo\\

Now, consider an irreducible pencil $\sigma\subseteq |D|$ on $X$.
For each integer $e$ such that
$0<e \leq \deg D$ let $\Lambda(\sigma,e)$ be the set of
pairs $W=({\mathcal L}, (v_p)_{p\in {BP}(\sigma)})$, where
${\mathcal L}$ is an effective class of $\pic(X)$ of degree $e$
and $(v_p)_{p\in {BP}(\sigma)}$ is a sequence of non-negative
integers, satisfying the following properties:

\begin{itemize}

\item[(a)] $v_p\leq e$ for all $p\in BP(\sigma)$,

\item[(b)] $v_p\geq \sum_q v_q$ for all $p\in BP(\sigma)$, where
the sum is taken over the set of points $q\in BP(\sigma)$ which belong to
the strict transform of the prime exceptional divisor associated with the blowing-up centered at
$p$,

\item[(c)] ${\mathcal L}^2\leq \sum_{p\in BP(\sigma)} v_p^2$,

\item[(d)] either $K_X\cdot {\mathcal L}+\sum_{p\in
BP(\sigma)}v_p\geq 0 $ and ${\mathcal L}^2+K_X\cdot {\mathcal L}+2
\geq \sum_{p\in BP(\sigma)} v_p(v_p-1)$, or ${\mathcal
L}^2=\sum_{p\in BP(\sigma)} v_p^2$ and $K_X\cdot {\mathcal
L}+\sum_{p\in BP(\sigma)}v_p=-2 $, or $K_X\cdot {\mathcal
L}+\sum_{p\in BP(\sigma)}v_p ={\mathcal L}^2- \sum_{p\in
BP(\sigma)} v_p^2=-1$, where $K_X$ stands for the canonical class
of $X$,

\item[(e)] $  {\mathcal L}\cdot D=\sum_{p\in {BP}(\sigma)} v_p
m_p(\sigma)$.

\end{itemize}

Taking into account the imposed conditions on the surface $X$ and
\cite[Lecture 16]{mum} one has that there exist finitely many
effective classes in $\pic(X)$ of degree $e$ and, therefore,
$\Lambda(\sigma,e)$ is finite.


\begin{alg}\label{alg}$ $\\
\noindent Input: $BP(\sigma)$.\\
\noindent Output: The set of integral components of the special
fibers of $\sigma$.

\noindent Begin
\begin{itemize}

\item[] Let $\Xi_0:= \emptyset$

\item[] Let $e:= 1$

\item[] While $e\leq d$ do:

\begin{itemize}

\item[] $\Xi_e:= \emptyset$

\item[] For each $W=({\mathcal L}, (v_p)_{p\in {BP}(\sigma)})\in \Lambda(\sigma,e)$ do:

\begin{itemize}

\item[]  If the complete linear system $|\pi_{\sigma}^*{\mathcal
L}-\sum_{p\in {BP}(\sigma)} v_p [E_p^*]|$ satisfies the following properties:

\begin{itemize}

\item[(1)]  it has projective dimension 0,

\item[(2)]  it has not exceptional part,

\item[(3)]  no curve in $\bigcup_{0 \leq j<e}\Xi_j$ is a component of the curve $C_W$ of $X$ given by the
direct image by $\pi_{\sigma}$ of the unique element of the
linear system,

\end{itemize}

then  let $\Xi_e:= \Xi_e\cup \{C_W\}$

\end{itemize}

\item[] Let $e:= e+1$

\end{itemize}

\item[] Return the set $\bigcup_{0< j\leq e}\Xi_j$

\end{itemize}

\noindent End

\end{alg}


Observe that, to check if a pair $W=({\mathcal L}, (v_p)_{p\in
{BP}(\sigma)})$ belonging to some set $\Lambda(\sigma,e)$ satisfies
the above properties $(1)$ and $(2)$ involves the knowledge of a
basis of the space of global sections of ${\mathcal L}$ and the
resolution of a system of linear equations. Also, a possibility for
checking the property $(3)$ is to verify that, for all $Q\in
\bigcup_{0\leq j<e}\Xi_j$, the linear system
$|\tilde{C}_W-\tilde{Q}|$ is empty (and, again, this can be done
solving a system of linear equations if one knows a basis of
$H^0(X,\co_X(C_W-Q))$.

The next result justifies the correctness of the algorithm.

\bl For each integer $e$ such that $1\leq e\leq \deg D$, $\Xi_e$ is the set of
integral components $C$ of special fibers of $\sigma$ such that $\deg C=e$.
\el

\noindent {\it Proof}. Assume that $C$ is an integral component of
degree $e$ of a special fiber of $\sigma$. Define ${\mathcal L}=[C]$
and $v_p:=m_p(C)$ for all $p\in {BP}(\sigma)$, and take
$W=({\mathcal L}, (v_p)_{p\in BP(\sigma)})$. Let us see that $W\in
\Lambda(\sigma,e)$. The first property defining this set, (a), must
be obviously satisfied by $W$, since the integers $v_p$ coincide
with the multiplicities of the strict transforms of the curve $C$.
(b) means that the proximity inequalities  \cite[Chap. II, book
4]{en}  are satisfied for $C$. (c) means that the self-intersection
of $\tilde{C}$ is non-positive, which is true by Proposition
\ref{gordo}. (d) follows from the Adjunction Formula and (e) follows
also from Proposition \ref{gordo}. We conclude, then, that $W\in
\Lambda(\sigma,e)$.

Taking into account again Proposition \ref{gordo}, each element of
the linear system $|\tilde{C}|=|{\pi_{\sigma}^*}{\mathcal
L}-\sum_{p\in {BP}(\sigma)} v_p [E_p^*]|$ is sum of strict
transforms of integral components of special fibers of the pencil
$\sigma$ and (possibly) strict transforms of exceptional divisors.
Hence, since the set of special fibers is finite, Property $(1)$ in
Algorithm \ref{alg} is satisfied for $W$. (2) holds clearly, since
$\tilde{C}$ belongs to the linear system. Therefore $C$ will be
obtained as $C_W$ in Algorithm \ref{alg}.

Conversely, set $C\in \Xi_e$ and let $W\in \Lambda(\sigma,e)$ be
such that $C=C_W$. Due to (e), Proposition \ref{gordo} and
Properties $(2)$ and $(3)$ in Algorithm \ref{alg}, it follows that
$C$ is an integral component of a fiber of the pencil $\sigma$. But,
taking into account $(1)$, it is clear that this
fiber must be special.\findemo\\

\noindent {\it Remark.} Notice that, from a basis of $\sigma$ and
the output of Algorithm \ref{alg}, one can determine which obtained
curves $C_W$ are components of the same special fiber and,
therefore, compute all the special fibers of the pencil.\\


\noindent {\it Example}. Consider the complex projective plane
$\gp^2_{\mathbb C}$, with projective coordinates $(X:Y:Z)$.
Let $H$ be a hyperplane section and consider the pencil
$\sigma\subseteq |3H|$  spanned by the divisors given by the
projective curves with equations
$F:=27X^3-27X^2Y+9XY^2-Y^3-8XZ^2+5YZ^2=0$ and
$G:=X^3+6X^2Y+12XY^2+8Y^3-7YZ^2=0$.
$BP(\sigma)$ consists of 9 points in $\gp^2_{\mathbb C}$ and
$m_p(\sigma)=1$ for all $p\in BP(\sigma)$. Therefore,
$G_{\sigma}=3H^*-\sum_{p\in BP(\sigma)} E_p^*$ and $\sigma$ is
irreducible (in other case, $\sigma$ would be composite with an
irreducible pencil, which contradicts the fact that $m_p(\sigma)=1$
for any point $p\in BP(\sigma)$). The conditions (a)-(e) given in
the above definition of the elements $W=({\mathcal L},(v_p)_{p\in
BP(\sigma)})$ of $\Lambda(\sigma,1)$ (resp., $\Lambda(\sigma,2))$
force $C_W$ to be a line (resp., conic) passing through three
(resp., six) points of $BP(\sigma)$, and $\Lambda(\sigma,3)$ is
empty. Using Algorithm \ref{alg} we can deduce that the integral
components of the special fibers of $\sigma$ are given by the
following equations:
$$L_1:=4X+Y=0;\;\;\; L_2:=2X-3Y=0;$$
$$L_3:=2X-(7\sqrt{5}+17)Y=0;\;\;\; L_4:=2X+(7\sqrt{5}-17)Y=0;$$ $$C_1:=(11\sqrt{5}-15)Y^2+2(2\sqrt{5}-25)XY+4(10+9\sqrt{5})X^2-16\sqrt{5}Z^2=0;$$
$$C_2:=(11\sqrt{5}+15)Y^2+2(2\sqrt{5}+25)XY+4(9\sqrt{5}-10)X^2-16\sqrt{5}Z^2=0;$$
$$C_3:=3Y^2+3XY+13X^2-4Z^2=0;\;\;\; C_4:=7Y^2-7XY+7X^2-2Z^2=0.$$
The special fibers of the pencil are those corresponding with the
equations:
$$F+G=L_1 C_4=0;\;\;\; F-G=L_2 C_3=0;$$
$$4\sqrt{5}F-4(9\sqrt{5}+20)G=L_3C_2=0;\;\;\;4\sqrt{5}F-4(9\sqrt{5}-20)G=L_4 C_1=0.$$

$ $\\

\end{document}